\documentclass[11pt]{article}
\usepackage{amsmath,amsfonts}
\newcommand{\bb}{\mathbb}

\newcommand{\cx}{{\bb C}}

\newcommand{\integers}{{\bb Z}}

\newcommand{\ratls}{{\bb Q}}
\newcommand{\reals}{{\bb R}}

\newcommand{\hthree}{{{\bb H}^3}}

\newcommand{\widemargins}{
\setlength{\textwidth}{5.8in}
\setlength{\oddsidemargin}{0.25in}
\setlength{\evensidemargin}{0.25in}
}

\newcommand{\qed}[1]{\nopagebreak[4]{\tiny \hfill
\fbox{\ref{#1}} \linebreak }\pagebreak[2]}
\newcommand{\del}{\partial}

\newcommand{\zbar}{{\overline{z}}}

\newcommand{\chat}{\widehat{\cx}}

\newcommand{\area}{\operatorname{area}}

\newcommand{\id}{\operatorname{id}}

\renewcommand{\Im}{\operatorname{Im}}

\newtheorem{theorem}{Theorem}[section]
\newtheorem{prop}[theorem]{Proposition}
\newtheorem{lemma}[theorem]{Lemma}
\newtheorem{cor}[theorem]{Corollary}
\newtheorem{conj}[theorem]{Conjecture}

\newcommand{\cA}{{\cal A}}
\newcommand{\cB}{{\cal B}}

\newcommand{\cH}{{\cal H}}

\newcommand{\cM}{{\cal M}}

\newcommand{\cO}{{\cal O}}
\newcommand{\cP}{{\cal P}}

\newcommand{\cR}{{\cal R}}
\newcommand{\cS}{{\cal S}}
\newcommand{\cT}{{\cal T}}

\widemargins
\renewcommand{\int}{{\operatorname{int}}}
\newcommand{\ab}{\cA\cB}

\newcommand{\pr}{{\operatorname{par}}}
\newcommand{\Ao}{{{\overset{\circ}{\cA_K}}}}
\newcommand{\Vo}{{{\overset{\circ}{V}}}}
\newcommand{\Vop}{{{\overset{\circ}{V'}}}}
\newcommand{\Uo}{{{\overset{\circ}{U}}}}
\newcommand{\Co}{{{\overset{\circ}{C}}}}
\newcommand{\Tu}{{\mathbb{T}}}

\begin{document}
\title{The space of Kleinian punctured torus groups is not locally connected}

\author{K. Bromberg\footnote{Supported by a grant from the NSF and a Sloan Foundation Fellowship}}

\date{\today}

\maketitle

\begin{abstract}
\noindent
We show that the space of Kleinian punctured torus groups is not locally connected.

\end{abstract}

\section{Introduction}
In this paper we give the first example of a deformation space of Kleinian groups that is not locally connected. Our space comes from one of the simplest examples of Kleinian groups, the punctured torus. These are Kleinian groups isomorphic to a two generator free group with the extra condition that the commutator is parabolic. Although there are some special features of punctured torus groups that play an important role in our proof the general method should be useful in studying other deformation spaces of Kleinian groups.


To put our result in perspective we give a brief historical survey of other results on the topology of deformation spaces of Kleinian groups. In our history we will restrict to results that describe unusual or exotic features of this topology.
\begin{itemize}
\item The Bers slice: There is a canonical bijection between any two Bers slices for a surface $S$. Kerckhoff and Thurston (\cite{Kerckhoff:Thurston}) showed that this bijection is in general not a homeomorphism and their methods lead one to believe that it is never a homeomorphism unless $S$ is a punctured torus or 4-times punctured sphere.

\item Bumping: For a general hyperbolizable 3-manifold $M$ the components of the interior of the deformation space will be in one to one correspondence with the marked homeomorphism types of $M$. Anderson and Canary (\cite{Anderson:Canary:pages}) found examples where these components bumped. In these examples distinct components of the interior of the deformation space have intersecting closure.
See also \cite{Anderson:Canary:McCulluogh:bump}, \cite{Holt:bump} and \cite{Holt:bigbump}.

\item Self-bumping: McMullen (\cite{McMullen:graft}) showed that quasifuchsian space self-bumps. That is there are points $\rho$ on the boundary such that the intersection of any sufficiently small neighborhood of $\rho$ with quasifuchsian space is disconnected. See also \cite{Ito:bump} and \cite{Bromberg:Holt:bump}.
\end{itemize}

The Sullivan dictionary encourages us to compare our result with what is known about rational maps of $\chat$. Following Sullivan, here is an abbreviated score card where we compare results on rational maps to results about Kleinian punctured torus groups. In our table $\hat{T}$ is the punctured torus.

\begin{table}[htdp]
\begin{center}\begin{tabular}{|c|c|c|}\hline   & \mbox{Kleinian groups} & \mbox{Rational maps of $\chat$} \\\hline \mbox{$1$-dim} & \mbox{Bers' slice of $\hat{T}$ and Maskit slice} & \mbox{Mandelbrot set} \\  \mbox{l. c.} & \mbox{Yes (\cite{Minsky:torus})} & \mbox{Main conjecture} \\\hline \mbox{$2$-dim} & $AH(\hat{T})$ & \mbox{Cubic connectivity locus} \\  \mbox{l. c.} & \mbox{No (this paper)} & \mbox{No (\cite{Lavaurs:cubics} and \cite{Milnor:cubics})} \\\hline \end{tabular} 
\end{center}
\label{defaulttable}
\end{table}
\vspace{-.22in}
\noindent The main point to notice is that in dimension one the deformation spaces are locally connected or at least conjectured to be. In dimension two the deformation spaces are not locally connected and we expect this to be true in all higher dimensions. Another point to emphasize is that the conjecture that the Mandelbrot set is locally connected is equivalent to giving a classification of quadratic polynomials (see \cite{Branner:mandelbrot}) analogous to Minsky's ending lamination theorem for punctured torus groups (\cite{Minsky:torus}). In the case of Kleinian groups we use this classification to prove non-local connectivity.

We close this introduction with some brief words on the proof. Our discussion will emphasize the importance of the 2nd dimension. Let $\cS$ be a closed subset of $\cx^2$ and assume that $\cS$ is invariant under the translation $(z,w) \mapsto (z, w+2)$. We add a collection of points at $\infty$ to $\cS$ by viewing $\cS$ as a subset of $\cx \times \chat$ and taking the closure. We then note that if for some $z$ the intersection of $\cS$ with $\{z\} \times \cx$ contains a bounded component then the closure of $\cS$ in $\cx \times \chat$ is not locally connected at $(z, \infty)$. Furthermore unless $\cS$ is very simple we can always find such a point $z$. We will prove our theorem by showing the deformation space of Kleinian punctured torus groups is locally homeomorphic to such a space at least at some cusp on the boundary of the deformations space. The translational symmetry will come from Dehn twisting along the curve corresponding to the cusp.

It is perhaps easier to visualize the phenomenon in $\reals^2$ rather than $\cx^2$. Again we assume that $\cS$ is a closed subset that is invariant under a vertical translation. We can compactify the vertical lines and take the closure of $\cS$ in $\reals \times S^1$. This set will not be locally connected unless the boundary of $\cS$ is vertical lines. On the other hand if we take $\cS$ to be a subset of $\cx$ and take the closure in $\chat$ then this closure will be locally connected as long as $\cS$ is locally connected and contains no bounded components.

{\bf Acknowledgements.} The main ideas for this paper were developed while the author was visiting the Newton Institute for the program on Spaces of Kleinian Groups and Hyperbolic 3-Manifolds. He would like to thank Caroline Series, Yair Minsky and Makoto Sakuma for organizing the program.

The author would also like to thank Jan and Nigel Evans for their hospitality at their home in Hahei where much of this paper was written.

The author also thanks Aaron Magid for his comments on a draft version of this paper.

\section{Preliminaries}
\subsection{Kleinian groups}
A {\em Kleinian group} is a discrete subgroup of $PSL_2\cx$. We are interested in the space of Kleinian groups isomorphic to a fixed group $G$. We begin by giving this space a topology.

Let $G$ be a finitely generated group and $\cR(G)$ the space of representations of $G$ in $PSL_2\cx$. If $G$ has $n$ generators then $\cR(G)$ embeds in $(PSL_2\cx)^n$ and this embedding gives $\cR(G)$ a topology. Let $R(G)$ be the quotient of $\cR(G)$ under the action of conjugation. In general the topology of $R(G)$ may be non-Hausdorff. However, as we'll see below, at the subspace that we are interested in $R(G)$ is a manifold.

Let $H(G)$ be the set of conjugacy classes of discrete, faithful representations of $G$ in $PSL_2\cx$. We give $H(G)$ the {\em algebraic topology} by taking the induced topology on $H(G)$ as a subset of $R(G)$. We denote $H(G)$ with this topology as $AH(G)$. As is standard practice we will often abuse notation and treat elements of $AH(G)$ as single representations instead of as conjugacy classes.

Let $N$ be a compact, hyperbolizable 3-manifold with boundary. If $G = \pi_1(N)$ then one denotes $R(G)$ and $AH(G)$ as $R(N)$ and $AH(N)$, respectively. Let $\cP$ be a collection of disjoint, essential and homotopically distinct annuli and all tori in $\del N$. Then $R(N, \cP) \subset R(N)$ and $AH(N, \cP) \subset AH(N)$ are those representations $\sigma$ where $\sigma(\gamma)$ is parabolic if $\gamma$ is freely homotopic into $\cP$. As a convention we will assume that all tori in $\del N$ are contained in $\cP$. In particular $\cP$ will be non-empty if $\del N$ contains tori. We also define $\del(N, \cP) = \del N - \cP$.

Here is the result, mentioned above, describing the topology of $R(G)$ at points in $AH(N, \cP)$:
\begin{theorem}[Kapovich \cite{Kapovich:book}]
\label{manifoldrepresentations}
A neighborhood of $AH(N, \cP)$ in $R(N, \cP)$ is a manifold whose dimension is the same as the dimension of the Teichm\"uller space of $\del(N, \cP)$.
\end{theorem}



\subsection{Hyperbolic geometry}
The Lie group $PSL_2\cx$ is isomorphic to orientation-preserving isometries of $\hthree$. Therefore the quotient $M_\sigma = \hthree/\sigma(\pi_1(N))$ is a complete hyperbolic 3-manifold. The representation also determines an isomorphism from $\pi_1(N)$ to $\pi_1(M_\sigma)$. We keep track of this isomorphism by picking a homotopy equivalence
$$f_\sigma: N \longrightarrow M_\sigma$$
such that $(f_\sigma)_*$ is the given isomorphism. The pair $(M_\sigma, f_\sigma)$ is a {\em marked, hyperbolic 3-manifold}. Note that the interior of $N$ need not be homeomorphic to $M_\sigma$. Even if it is homeomorphic the map $f_\sigma$ may not be a homeomorphism or even homotopic to a homeomorphism.

Conversely a marked, hyperbolic 3-manifold $(M, f)$ determines a representation $\sigma$. Namely the hyperbolic manifold $M$ can be written as a quotient $M = \hthree/\Gamma$ where $\Gamma$ is a group of isometries of $\hthree$ isomorphic to $\pi_1(N)$. The marking map $f$ determines an isomorphism between $\pi_1(N)$ and $\pi_1(M)$ via the map $f_*$. This gives an isomorphism between $\pi_1(N)$ and $\Gamma$ which is our representation. We can define an equivalence relation for marked, hyperbolic 3-manifolds by setting $(M_0, f_0) \sim (M_1, f_1)$ if there is an orientation preserving isometry $g: M_0 \longrightarrow M_1$ such that $g \circ f_0$ is homotopic to $f_1$. Two marked, hyperbolic 3-manifolds will be equivalent if and only if the induced representations are conjugate. There is therefore a bijection between equivalence classes of marked, hyperbolic 3-manifolds and conjugacy classes of representations in $AH(N)$. Throughout the paper we will implicitly use this bijection and refer to points in $AH(N)$ as both representations and marked, hyperbolic 3-manifolds.

\subsection{The thick-thin decomposition}

Given $\epsilon>0$ let $M^{\leq \epsilon}$ be the set of points in the hyperbolic manifold $M$ whose injectivity radius is $\leq \epsilon$. We can similarly define $M^{\geq \epsilon}$. By the Margulis lemma (\cite{Kazhdan:Margulis}) there is a constant $\epsilon_M$ such that for any $\epsilon \leq \epsilon_M$ the components of $M^{\leq \epsilon}$ have three possible types:
\begin{enumerate}
\item a tubular neighborhood of a geodesic of length $< 2\epsilon$ (or a closed geodesic of length $2\epsilon$),

\item the quotient of a horoball by a rank one parabolic subgroup or

\item the quotient of a horoball by a rank two parabolic subgroup.
\end{enumerate}
The last two types are the {\em parabolic} thin part of $M$ which we denote $M^{\pr, \leq \epsilon}$.

If $\gamma$ is a simple closed curve of length $<\epsilon \leq \epsilon_M$ then $\Tu^{\epsilon}(\gamma,M)$ is the component of $M^{\leq \epsilon}$ that contains $\gamma$. This is case (1) and $\Tu^{\epsilon}(\gamma, M)$ is a {\em Margulis tube}. Thin parts of types (1) and (2) are associated to components of the paring locus. If $A$ is a annulus in the paring locus of $M$ then there is a corresponding component of $M^{\leq \epsilon}$ of type (2) and it is denoted $\Tu^{\epsilon}(A, M)$. If $T$ is a torus in the paring locus then the corresponding component of $M^{\leq \epsilon}$ is of type (3) and is denoted $\Tu^{\epsilon}(T,M)$. In case (2), $\Tu^{\epsilon}(A, M)$ is a {\em rank one cusp}. In case (3), $\Tu^{\epsilon}(T,M)$ is a {\em rank two cusp}.

Our final goal of this subsection is to give a description of the geometry of a rank two cusp $\Tu^{\epsilon}(T,M)$. We start by noting that a horoball is foliated by horospheres and in the quotient cusp these horospheres descend to Euclidean tori.  Although these tori are not all isometric they are all conformally equivalent. Furthermore the conformal structure of this tori completely determine the cusp. Our purpose will be to use this conformal structure to define a notion of length for curves on $T$.


Let $\beta$ be a simple closed curve on $T$ and choose a Euclidean torus $T'$ in the foliation of $\Tu^{\epsilon}(M,T)$ described in the previous paragraph. There is a canonical homotopy equivalence between $T$ and $T'$ and $\beta$ will be sent to a simple closed curve $\beta'$ under this homotopy equivalence. We then define the normalized length of $\beta$ to be
$$\frac{\operatorname{length}(\beta')}{\sqrt{\area(T')}}$$
where $\operatorname{length}(\beta')$ is the length of a geodesic representative of $\beta'$. This normalized length does not depend on the choice of torus $T'$ and only on the rank two cusp $\Tu^{\epsilon}(M, T)$.

To calculate this length we can assume that the rank two cusp has generators
$$\left(\begin{array}{cc} 1 & 2 \\ 0 & 1 \end{array} \right) \mbox{ and }
\left(\begin{array}{cc} 1 & z \\ 0 & 1 \end{array} \right)$$
where $\Im z > 0$. If the curve $\beta$ corresponds to the second generator then its normalized length is
$$\frac{|z|}{\sqrt{2\Im z}}.$$

\subsection{The Ahlfors-Bers parameterization}
\label{abparam}
The Ahlfors-Bers parameterization is a parameterization of the interior of $AH(N, \cP)$ by certain Teichm\"uller spaces. This a classical result and a standard reference is \cite{Bers:def} or \cite{Canary:McCullough:book}.

Although the parametization is defined for any $N$ with with finitely generated fundamental group we will restrict ourselves here to manifolds with incompressible boundary. That is we will assume that for each component $S$ of $\del(N,\cP)$ the inclusion of $S$ in $N$ is injective on the level of fundamental groups.

To define this parameterization we also need some standard notions from Kleinian groups. Let $\sigma$ be a representation in $AH(N, \cP)$. Then the action of $\sigma$ on $\hthree$ extends continuously to an action on the boundary of $\hthree$ which can be naturally identified with the Riemann sphere, $\chat$. The {\em domain of discontinuity}, $\Omega_\sigma$, is the largest subset of $\chat$ where $\sigma(\pi_1(N))$ acts properly discontinuously. Then the quotient $\bar{M}_\sigma = (\hthree \cup \Omega_\sigma)/\sigma(\pi_1(N))$ is a 3-manifold with boundary. The interior of $\bar{M}_\sigma$ is the original hyperbolic manifold $M_\sigma$. By Ahlfors' Finiteness Theorem, the boundary of $\bar{M}_\sigma$ is a finite collection of finite area hyperbolic Riemann surfaces. This is the {\em conformal boundary} of $M_\sigma$.

The {\em convex core}, $C(M_\sigma)$, of $M_\sigma$ is the smallest convex subset of $M_\sigma$ where the inclusion map is a homotopy equivalence. The manifold $M_\sigma$ is {\em geometrically finite} if $C(M_\sigma)$ has finite volume. Let $GF(N, \cP)$ be those representations $\sigma$ in $AH(N, \cP)$ where $M_\sigma$ is geometrically finite.

A representation $\sigma \in AH(N, \cP)$ is {\em minimally parabolic} if $\sigma(g)$ is parabolic if and only if $g$ is freely homotopic into $\cP$. It is convenient to let $MP(N, \cP) \subset AH(N, \cP)$ be the set of representations that are both minimally parabolic and geometrically finite.

With these definitions in place we can now describe the Ahlfors-Bers parameterization. The first result we will need is the following:
\begin{theorem}[Marden \cite{Marden:kgs}, Sullivan \cite{Sullivan:QCDII}]
\label{interior}
The interior of $AH(N, \cP)$ (as a subspace of $R(N, \cP)$) is $MP(N, \cP)$.
\end{theorem}

In general $MP(N, \cP)$ will have many components and possibly even an infinite number.  We will use the marking maps to pick out a distinguished component of $MP(N, \cP)$.
First choose an orientation for $N$ and let $AH_0(N, \cP) \subseteq AH(N, \cP)$ be the set of representations $\sigma$ where the marking map $f_\sigma$ can be chosen to be an orientation preserving embedding and $f_\sigma(N) \cap M^{\pr, \leq \epsilon} = f_\sigma(\cP)$ for some $\epsilon < \epsilon_M$. We then define $MP_0(N, \cP)$ to be the intersection of $AH_0(N, \cP)$ and $MP(N, \cP)$.

We first define the parameterization for $MP_0(N, \cP)$. The key point is that if $\sigma$ is in $MP_0(N,\cP)$ then the marking map $f_\sigma$ can be homotoped to a homeomorphism from $N - \cP$ to $\bar{M}_\sigma$. If we restrict this map to a component $S$ of $\del(N, \cP)$ then the pair $(S, f_\sigma)$ determines a point in the Teichm\"uller space of $S$. If the boundary is incompressible this point in Teichm\"uller space does not depend on our choice of marked, hyperbolic 3-manifold in the equivalence class correspoding to $\sigma$. Note that if $S$ is a triply-punctured sphere then the Teichm\"uller space of $S$ is a point so we will ignore such components.
By repeating this construction for all other components we define an {\em Ahlfors-Bers parameterization}
$$\ab : MP_0(N, \cP) \longrightarrow \cT(\del(N, \cP))$$
where $\cT(\del(N, \cP))$ is the product of the Teichm\"uller spaces of the components of $\del(N, \cP)$ that are not triply-punctured spheres.

We would like to extend this parameterization to the rest of the components of $MP(N, \cP)$. This can be done as follows.  Let $\sigma$ be a representation in a component $B$ of $MP(N, \cP)$ and let $(N', \cP')$ be a pared 3-manifold such that $N' - \cP'$ is homeomorphic to $\bar{M}_\sigma$. We also choose a homeotopy equivalence
$f'_\sigma : N' \longrightarrow M_\sigma$
that can be homotoped to an embedding as in the definition of $AH_0$. In particular the pair $(M_\sigma, f'_\sigma)$ induces a representation $\sigma'$ in $AH(N', \cP')$ and this representation lies in $MP_0(N', \cP')$. Note that the marking maps induce isomorphisms on the fundamental group so the composition $(f_\sigma)^{-1}_* \circ (f'_\sigma)_*$ is an isomorphism from $\pi_1(N')$ to $\pi_1(N)$. This defines an homeomorphism from $AH(N, \cP)$ to be $AH(N', \cP')$ by pre-composing the representations in $AH(N, \cP)$ with $(f_\sigma)^{-1}_* \circ (f'_\sigma)_*$. This homeomorphism will take $\sigma$ to $\sigma'$ and therefore will take the component $B$ of $MP(N, \cP)$ to $MP_0(N', \cP')$. We then use the Ahlfors-Bers parameterization for $MP_0(N', \cP')$ to define the Ahlfors-Bers parameterization for $B$.





By repeating this construction for all components we can define $\ab$ on all of $MP(N, \cP)$. The image of $\ab$ will then be a disjoint union of products of Teichm\"uller spaces. We have the following foundational result.
\begin{theorem}[Ahlfors, Bers, Kra, Maskit, Marden, Sullivan]
\label{param}
The map $\ab$ is a homeomorphism.
\end{theorem}
For manifolds with compressible boundary a similar statement holds where the image of the Ahlfors-Bers parameterization is a certain quotient of Teichm\"uller space. We will not address manifolds with compressible boundary in this paper.

\subsection{Drilling and filling}


One of our main tools will be the drilling and filling theorems. These two results compare hyperbolic metrics on a manifold and on the same manifold with a closed geodesic removed. Their proofs involve the deformation theory for hyperbolic cone-manifolds developed in \cite{Hodgson:Kerckhoff:cone} and \cite{Hodgson:Kerckhoff:dehn} and extended to geometrically finite manifolds in \cite{Bromberg:rigidity} and \cite{Bromberg:schwarz}. Before stating the theorems we need to setup some notation.

Let $N$ be a compact 3-manifold with boundary, $\gamma$ a simple closed curve in the interior of $N$ and $W$ an open tubular neighborhood of $\gamma$. Let $\hat{N} = N - W$ and let $T = \del W$ be the torus in the boundary of $\hat{N}$. Let $\beta$ be an essential simple closed curve on $T$ that bounds a disk in $W$. The curve $\beta$ is unique up to homotopy. We would like to compare hyperbolic metrics on $N$ with hyperbolic metrics on $\hat{N}$. In general, metrics on these manifolds will not be unique so to get an interesting comparison we need to be careful when choosing them. The conformal boundary gives us a way of doing this.

Let $M$ be a hyperbolic structure on the interior of $N$ such that $\gamma$ is a closed geodesic. Let $\hat{M}$ be a hyperbolic structure on the interior of $\hat{N}$. Let $\phi: \hat{M} \longrightarrow M$ be the inclusion map. Assume that the hyperbolic structure $\hat{M}$ has been chosen such that $\phi$ extends to a conformal map between the conformal boundaries of $M$ and $\hat{M}$. 
%
We say that $\hat{M}$ is the {\em $\gamma$-drilling} of $M$ and $M$ is the {\em $\beta$-filling} of $\hat{M}$. Note that if we are given $\hat{M}$ there may not exist a $\beta$-filling $M$. On the other hand for any hyperbolic 3-manifold $M$ and simple closed geodesic $\gamma$ Kerckhoff has observed that there is a $\gamma$-drilling $\hat{M}$ (see \cite{Kojima:cone}).

In the following result we obtain bounds on the bi-Lipschitz constant of of the inclusion map $\phi$.
\begin{theorem}
\label{drilling}
Given any $L>1$, $K>0$  and $\epsilon>0$ there exists an $\ell>0$ such that the following holds. Let $M$ be a hyperbolic 3-manifold and $\gamma$ a simple closed geodesic in $M$ with length $\leq \ell$. Let $\hat{M}$ be the $\gamma$-drilling of $M$.
\begin{enumerate}
\item The inclusion map $\phi$ can be chosen such that $\phi$ restricts to an $L$-bi-Lipschitz diffeomorphism from $M - \Tu^{\epsilon}(\gamma, M)$ to $\hat{M} - \Tu^{\epsilon}(T, \hat{M})$.

\item Let $\beta$ be a essential simple closed curve  in the new cusp of $\hat{M}$ that bounds a disk in $M$.  The normalized length of $\beta$ is $> K$.
\end{enumerate}
\end{theorem}

{\bf Remark.} In \cite{Brock:Bromberg:density} only (1) of Theorem \ref{drilling} was proven and only for manifolds without rank one cusps. This result can be extended to manifolds with rank one cusps using Theorem 3.4 in \cite{Brock:Bromberg:Evans:Souto}. Property (2) can be derived from the estimates on pp. 406-7 of \cite{Hodgson:Kerckhoff:dehn}.

\medskip

Next we state the filling theorem. Here we want the normalized length of $\beta$ to be long to get a small bi-Lipschitz constant.
\begin{theorem}
\label{filling}
Given any $L>1$, $\ell>0$ and $\epsilon>0$ there exists an $K>0$ such that the following holds. Let $\hat{M}$ be a hyperbolic 3-manifold with a rank two cusp and let $\beta$ be an essential simple closed in the cusp with normalized length $\geq K$.
\begin{enumerate}
\setcounter{enumi}{-1}
\item The $\beta$-filling $M$ of $\hat{M}$ exists. 

\item Let $\gamma$ be the the geodesic representative of the core curve of the filled torus in $M$. Then the filling map $\psi$ can be chosen such that it restricts to an $L$-bi-Lipschitz diffeomorphism from $\hat{M} - \Tu^{\epsilon}(T, \hat{M})$ to $M - \Tu^{\epsilon}(\gamma, M)$.

\item The length of $\gamma$ in $M$ is $< \ell$.
\end{enumerate}
\end{theorem}

{\bf Remark.} In \cite{Hodgson:Kerckhoff:dehn} (0) and (2) of Theorem \ref{filling} were proven for finite volume manifolds. The extension of this result to geometrically finite manifolds without rank one cusps follows directly from the methods of \cite{Bromberg:rigidity} and \cite{Bromberg:schwarz}. We can add rank one cusps exactly as for the drilling theorem by using Theorem 3.4 of \cite{Brock:Bromberg:Evans:Souto}. Finally we note that (1) of Theorem \ref{drilling} and (0) and (2) of Theorem \ref{filling} together imply (1) of Theorem \ref{filling}.

\subsection{Hyperbolic metrics on compact manifolds}

In this section we study hyperbolic metrics on a compact manifold $N$ with boundary. Let $\cH(N)$ be the space of smooth, complete hyperbolic metrics on $N$. Here by smooth we mean that the metric extends to a smooth metric on the boundary of $N$. We give $\cH(N)$ the $C^\infty$-topology. A basis for this topology is sets of the following form: Fix a metric $g$ and a constant $L>1$. Then a set in the basis is formed by taking all metrics $g'$ in $\cH(N)$ where the identity map is an $L'$-bi-Lipschitz map from the $g$-metric to the $g'$-metric for some  $L' \in [1, L)$.

A hyperbolic metric on $N$ determines a holonomy representation of $\pi_1(N)$ which is defined up to conjugation. Let
$$h : \cH(N) \longrightarrow R(N)$$
be the map which assigns to each hyperbolic metric its holonomy representation. Note that $\cH(N)$ will be an infinite dimensional space and that the map will not be injective. The following result gives a description of a neighborhood of a metric $g \in \cH(N)$ in terms of a neighborhood of the holonomy representation.
\begin{theorem}[\cite{Canary:Epstein:Green}]
\label{holonomyinverse}
Let $N_0 \subset N$ be the complement of an open collar neighborhood of $N$. Let $g$ be a hyperbolic metric in $\cH(N)$ with holonomy representation $\rho = h(g)$. Then there are neighborhoods $V$ of $\rho$ and $U$ of $g$ with the following properties:
\begin{enumerate}
\item On $V$ there exists a holonomy inverse
$$h^{-1} : V \longrightarrow \cH(N)$$
such that $h \circ h^{-1} = \id$ and $h^{-1}(\rho) = g$.

\item For $U$ we have $h(U) \subset V$ and for each $g' \in U$ there exists a diffeomorphism
$$\alpha : N \longrightarrow N$$
isotopic to the identity sucht that $h^{-1} \circ h(g')$ and $\alpha_* g'$ agree on $N_0$.
\end{enumerate}
\end{theorem}

\section{The model space}
We will study the topology of $AH(N, \cP)$ at geometrically finite boundary points that contain one extra conjugacy class of parabolics. We will only work with manifolds with incompressible boundary and if the reader likes one can assume that $N = S \times [0,1]$ and $\cP = \del S \times [0,1]$ where $S$ is a compact surface with boundary. However, as it is no extra work, we will allow general pared manifolds with incompressible boundary.  We begin by setting the notation that we will use for the remainder of this section.

Let $(N, \cP)$ be a pared, compact, hyperbolizable 3-manifold with pared incompressible boundary. Let $\gamma$ be a simple closed curve on $\del (N, \cP)$ that is not homotopic into $\cP$. Let $b$ be an element of $\pi_1(N)$ that is freely homotopic to $\gamma$ and assume that $b$ is primitive. Let $W$ be a open solid torus in the interior of $N$ whose core curve is isotopic to $\gamma$ and let $T$ be the boundary of $W$. Then define $\hat{N} = N \backslash W.$ Then $\pi_1(\hat{N})$ has presentation
%

$$\langle \pi_1(N), c | [b,c] \rangle.$$

Let $\cP'$ be the union of $\cP$ and an annulus with core curve $\gamma$. Let $\hat{P}$ be the union of the image of $\cP$ under the inclusion of $N$ in $\hat{N}$ and the torus $T$. Note that every representation in $MP_0(N, \cP')$ is also a representation in $AH(N, \cP)$. Our goal is to give a model for the topology of $AH(N, \cP)$ in a neighborhood of a representation in $MP_0(N, \cP')$.

Let $\sigma$ be a representation of $\pi_1(N)$ such that
$$\sigma(b) = \left( \begin{array}{cc}1 & 2 \\ 0 & 1 \end{array} \right).$$
To extend $\sigma$ to a representation of $\pi_1(\hat{N})$  of we only need to define the extension on $c$ and we need this new element to commute with $\sigma(b)$. There is a one complex dimension family of elements of $PSL_2\cx$ that will commute with $\sigma(b)$ and we can therefore define a one parameter family of extensions of $\sigma$. For each $z \in \cx$ we extend $\sigma$ to a representation $\sigma_z$ of $\pi_1(\hat{N})$ by setting
$$\sigma_z(c) = \left( \begin{array}{cc} 1 & z \\ 0 & 1 \end{array} \right).$$

\begin{lemma}
\label{standardconjugation}
Every representation $\sigma$ in $\cR(N, \cP')$ is conjugate to a representation $\sigma'$ with
$$\sigma'(b) = \left( \begin{array}{cc}1 & 2 \\ 0 & 1 \end{array} \right).$$
Further if $\sigma$ and $\sigma'$ are conjugate representations with
$$\sigma(b) = \sigma'(b) = \left( \begin{array}{cc}1 & 2 \\ 0 & 1 \end{array} \right)$$
then $\sigma_z$ and $\sigma'_z$ are also conjugate.
\end{lemma}

{\bf Proof.} The first statement follows from the fact that all parabolics in $PSL_2\cx$ are conjugate.

For the second statement we note that any conjugation which fixes the parabolic $\sigma(b)$ will also fix everything that commutes with $\sigma(b)$. Therefore the same element of $PSL_2\cx$ that conjugates $\sigma$ to $\sigma'$ will also conjugate $\sigma_z$ to $\sigma'_z$. \qed{standardconjugation}

We now define a map
$$\tau: MP_0(N, \cP') \times \cx \longrightarrow R(\hat{N}, \hat{\cP})$$
by $\tau(\sigma, z) = \sigma_z$. By Lemma \ref{standardconjugation} this map is well defined.
\begin{lemma}
\label{taucontinuity}
The map $\tau$ is a local homeomorphism at all points $(\sigma, z)$ such that $\sigma_z \in AH(\hat{N}, \hat{\cP})$.
\end{lemma}

{\bf Proof.} The map is easily seem to be continuous and injective. The statement then follows from Theorem \ref{manifoldrepresentations} and invariance of domain. \qed{taucontinuity}

Next, for each $\sigma \in MP_0(N, \cP')$ we define the set $\cA_{\sigma} \subset \cx$ by
$$\cA_{\sigma} = \{z \in \cx |  \tau(\sigma, z) \in AH(\hat{N}, \hat{\cP}) \mbox{ and } \Im z \geq 0\}.$$
Define $\cA \subset MP_0(N, \cP') \times \chat$ by
$$\cA = \{(\sigma, z) | z \in \cA_{\sigma} \mbox{ or $z = \infty$}\}.$$

We will construct a local homeomorphism from a neighborhood of $(\sigma, \infty)$ in $\cA$ to $AH(N, \cP)$. On the interior of $\cA$ this construction will work in a general setting but to extend the map to the boundary of $\cA$ we will need to restrict to the case of the punctured torus.



\subsection{From $\cA$ to $AH(N, \cP)$}
We begin by defining a map $\Phi$ from $\cA$ to $AH(N, \cP)$.
Recall from the filling theorem that there exists a $K>0$ such that any curve in a rank two cusp of normalized length $\geq K$ can be filled. Let $\cA_K \subset \cA$ be those pairs $(\sigma, z) \in \cA$ where
$$\frac{|z|}{\sqrt{2 \Im z}} > K.$$
Note that as $z \rightarrow \infty$ in $\chat$ the quantity on the left also limits to $\infty$. It is therefore natural to assume that $(\sigma, \infty)$ is in $\cA_K$ for all $\sigma \in MP_0(N, \cP')$. Ultimately we would like the domain of definition of $\Phi$ to be all of $\cA_K$. In this section we will only define $\Phi$ at those points $(\sigma, z)$ where either $\sigma_z \in MP(\hat{N}, \hat{\cP})$ or $z = \infty$. Note that the set of such points is dense in $\cA_K$ and that it contains the interior of $\cA_K$. We denote this set by $\Ao$.

If $(\sigma, z) \in \Ao$ let $\hat{M}_{\sigma,z}$ be the quotient hyperbolic manifold for the representation $\sigma_z$. Note that the normalized length of $c$ is $> K$ so $\hat{M}_{\sigma,z}$ will be $c$-fillable. Let $M_{\sigma,z}$ be the $c$-filling and let
$$\phi_{\sigma,z} : \hat{M}_{\sigma,z} - \Tu^{\epsilon_M}(T, \hat{M}_{\sigma,z}) \longrightarrow M_{\sigma, z} - \Tu^{\epsilon_M}(\gamma, M_{\sigma,z})$$
be the filling map. We also have a marking map
$$f_\sigma : N \longrightarrow M_\sigma$$
where $M_\sigma$ is the quotient hyperbolic manifold for the representation $\sigma$. As $\sigma_z$ is an extension of $\sigma$ the hyperbolic manifold $M_\sigma$ covers $\hat{M}_{\sigma, z}$. Let
$$\pi_{\sigma, z} : M_\sigma \longrightarrow \hat{M}_{\sigma, z}$$
be the covering map. Let $f_{\sigma, z} = \phi_{\sigma,z} \circ \pi_{\sigma, z} \circ f_\sigma$. 

%




\begin{lemma}
\label{newmarkingmap}
The map $(f_{\sigma, z})_*$ is an isomorphism and therefore $(M_{\sigma, z}, f_{\sigma, z})$ is a marked hyperbolic 3-manifold in $AH(N, \cP)$.
\end{lemma}

{\bf Proof.} The fundamental group $\pi_1(\hat{M}_{\sigma,z})$ has presentation
$$\langle \pi_1(M_\sigma), c | [b,c] \rangle$$
and the image of $\pi_1(N)$ under the composition $\pi_{\sigma,z} \circ f_\sigma$ is $\pi_1(M_\sigma)$. By the Seifert-Van Kampen Theorem the filling map $(\phi_{\sigma,z})_*$ is the quotient map
$$\langle \pi_1(M_{\sigma}), c | [b,c] \rangle \longrightarrow \langle \pi_1(M_{\sigma}), c | [b,c] \rangle/\langle c\rangle$$
where $\langle c \rangle$ is the normal subgroup of $\pi_1(M_{\sigma, z})$ generated by $c$. This quotient is exactly $\pi_1(M_\sigma)$ so $(f_{\sigma,z})_*$ is an isomorphism as desired. \qed{newmarkingmap}

We can now define $\Phi$ by
$$\Phi(\sigma, z) = \left\{ \begin{array}{cl} (M_{\sigma,z}, f_{\sigma,z}) & \mbox{if $z \neq \infty$} \\ \sigma & \mbox{if $z = \infty$.} \end{array} \right.$$




\subsection{From $AH(N, \cP)$ to $\cA$}
We now construct a map $\Psi$ from $AH(N, \cP)$ to $\cA$. To do so we fix a representation $\rho$ in $MP_0(N, \cP')$ and our map will be defined only in a neighborhood of $\rho$. We will need to choose this neighborhood carefully and although we will define a neighborhood of $\rho$ in $AH(N, \cP)$ in this section we will only be able to define $\Psi$ on the restriction of this neighborhood to $MP_0(N, \cP') \cup MP(N, \cP)$. To extend $\Psi$ to the entire neighborhood we will have to restrict to the special case of the punctured torus.

We now describe the neighborhood of definition of $\Psi$. To do this we first fix a smooth embedding
$$s_\rho : N \longrightarrow M_\rho$$
and let $g_\rho$ be the pullback via $s_\rho$ of the hyperbolic metric on $M_\rho$. We choose $s_\rho$ such that each annular component of $\cP'$ has core curve with length $< \epsilon_M/4$ in the $g_\rho$-metric. In particular, $\gamma$ has length $< \epsilon_M/4$ in the $g_\rho$-metric. With this is in place we can now define a neighborhood $V'$ of $\rho$ in $AH(N,\cP)$ with the following properties:
\begin{enumerate}
\item There exists an inverse
$$h^{-1} : V' \longrightarrow \cH(N)$$
to the holonomy map such that $h^{-1}(\rho) = g_\rho$. The existence of a neighborhood with such an inverse is given by (1) of Theorem \ref{holonomyinverse}.

\item The identity map on $N$ is a $2$-bi-Lipschitz map between any two metrics in $h^{-1}(V')$.

\item Let $g_\sigma = h^{-1}(\sigma)$ for $\sigma \in V'$ and let
$$s_\sigma: N \longrightarrow M_{\sigma}$$
be a local isometry from the $g_\sigma$-metric on $N$ to the hyperbolic metric on $M_\sigma$. Then there exists an $\epsilon_0 < \epsilon_M$ such that $s_\sigma(N) \subset M^{\geq \epsilon_0}_\sigma$.

\item Let $\ell$ be the constant given by the drilling theorem (Theorem \ref{drilling}) such that for any geodesic of length $< \ell$  the drilling map is uniformly bi-Lipschitz outside of the $\epsilon_0$-Margulis tube about the geodesic. For any $\sigma \in V'$, the length of $\gamma$ in $M_\sigma$ is less then $\min\{\epsilon_{M}/8, \ell\}$. This is possible since the length of $\gamma$ is zero in $M_\rho$ and the length function is continuous on $AH(N, \cP)$.

\item The intersection $V' \cap AH(N, \cP')$ is contained in $MP_0(N, \cP')$. This is possible since $MP_0(N, \cP')$ is open in $AH(N, \cP')$.

\end{enumerate}

Let $\Vop = V' \cap (MP_0(N, \cP') \cup MP(N, \cP))$.
For each $\sigma \in \Vop$, $\Psi(\sigma)$ will be a pair consisting of a representation and a complex number. To define the first half of $\Psi$ we construct a map
$$\omega: \Vop \longrightarrow AH(N, \cP').$$
As with the construction of $\Phi$ there will be two cases depending on whether or not $\sigma(b)$ is parabolic. 

We first assume that $\sigma(b)$ is not parabolic. This is the more involved case. Let $\gamma_\sigma$ be the geodesic representative of $\gamma$ in $M_\sigma$.  Let $\hat{M}_\sigma$ be the $\gamma_\sigma$-drilling of $M_\sigma$. Let
$$\psi_\sigma : M_\sigma - \Tu^{\epsilon_0}(\gamma_\sigma, M_\sigma) \longrightarrow \hat{M}_\sigma - \Tu^{\epsilon_0}(T, \hat{M}_\sigma)$$
be the drilling map. Let $\bar{M}_\sigma$ be the cover of $\hat{M}_\sigma$ induced by the image of $\pi_1(N)$ under the map $(\psi_\sigma \circ f_\sigma)_*$. Let
$$\bar{f}_\sigma: N \longrightarrow \bar{M}_\sigma$$
be the lift of $\psi_\sigma \circ f_\sigma$. 

\begin{lemma}
\label{parabolics}
The pair $(\bar{M}_\sigma, \bar{f}_\sigma)$ is a marked hyperbolic 3-manifold $(\bar{M}_\sigma, \bar{f}_\sigma)$ determining a representation in $AH(N, \cP')$.
\end{lemma}

{\bf Proof.} The representation induced by $(\bar{M}_\sigma, \bar{f}_\sigma)$ is clearly discrete since its image lies in a Kleinian group. Faithfulness comes from that fact that $(\psi_{\sigma} \circ f_\sigma)_*$ is injective. 

We now show that we have the correct cusps. That is we need to show that $\bar{f}_\sigma$, or equivalently $\psi_\sigma \circ f_\sigma$, can be homotoped to take $\cP'$ into cusps. By the properties of the filling map we can homotop $\psi_\sigma \circ f_\sigma$ to take $\cP'$ into cusps if we can homotop $f_\sigma$ in $M_\sigma - \Tu^{\epsilon_0}(\gamma_\sigma, M_\sigma)$ to take each component of $\cP'$ into a cusp of $M_\sigma$ or into $\Tu^{\epsilon_M}(\gamma_\sigma, M_\sigma) - \Tu^{\epsilon_0}(\gamma_\sigma, M_\sigma)$.

We only need to worry about the annular components of $\cP'$ since the tori will automatically be homotopic into cusps. Each annulus in $\cP'$ will have core curve of length $< \epsilon_M/2$ in the $g_\sigma$-metric by our choice of $f_\rho$ and (2). Therefore these core curves will be mapped into $M^{\leq \epsilon_M}_\sigma$. Those annuli that also lie in $\cP$ must have core curves mapped into cusps since $\sigma$ is in $AH(N, \cP)$. The one remaining annulus in $\cP' - \cP$ has core curve $\gamma$. It will be mapped into $\Tu^{\epsilon_M}(\gamma_\sigma, M_\sigma)$ and since the image of $f_\sigma$ is contained in $M_{\sigma}^{\geq \epsilon_0}$ we have $f_\sigma(\gamma) \subset \Tu^{\epsilon_M}(\gamma_\sigma, M_\sigma) - \Tu^{\epsilon_o}(\gamma_\sigma, M_\sigma)$. Therefore the representation induced by $(\bar{M}_\sigma, \bar{f}_\sigma)$ is in $AH(N, \cP')$. \qed{parabolics}



We can now define $\omega$ by
$$\omega(\sigma) = \left\{ \begin{array}{cl} (\bar{M}_\sigma, \bar{f}_\sigma) & \mbox{if $\sigma(b)$ is not parabolic} \\ \sigma & \mbox{if $\sigma(b)$ is parabolic.} \end{array} \right.$$

\begin{lemma}
\label{barcontinuous}
The map $\omega$ is continuous at all points in $\Vop \cap MP_0(N, \cP')$.
\end{lemma}

{\bf Proof.} Let $\sigma_i$ be a sequence in $\Vop$ that converges to $\sigma \in MP_0(N, \cP')$.
We will show that the sequence $\omega(\sigma_i)$ converges to $\omega(\sigma) = \sigma$. We can assume that the $\sigma_i$ are not in $\Vop \cap MP_0(N, \cP')$ because on $\Vop \cap MP_0(N, \cP')$ the map is the identity and therefore continuous. 

Note that the metrics $g_{\sigma_i}$ converge to the metric $g_\sigma$ in $\cH(N)$ since $h^{-1}$ is continuous. Let $\bar{g}_{{\sigma}_i}$ be the pullback via $\bar{f}_{\sigma_i}$ of the hyperbolic metric on $\bar{M}_{\sigma_i}$. As $\sigma_i$ limits to $\sigma$ the length of $\gamma_{\sigma_i}$ limits to zero. Therefore the bi-Lipschitz constants of the drilling maps $\psi_{\sigma_i}$ will limit to $1$. This implies
$$\underset{i \rightarrow \infty}{\lim}\bar{g}_{{\sigma}_i} = \underset{i \rightarrow \infty}{\lim} g_{\sigma_i} = g_\sigma$$
and therefore we must have $\omega(\sigma_i) \rightarrow \sigma$ as desired. \qed{barcontinuous}

We can now finish the definition of $V$, the neighborhood of $\rho$ where $\Psi$ will be defined.
\begin{cor}
\label{minimalparabolics}
There exists a neighborhood $V \subset V'$ of $\rho$ such that if $\sigma \in \Vo = V \cap (MP_0(N, \cP') \cup MP(N, \cP))$ then $\omega(\sigma) \in MP_0(N, \cP')$.
\end{cor}

{\bf Proof.} Assume not. Then every neighborhood of $\rho$ in $AH(N, \cP)$ contains a representation $\sigma$ in $MP(N, \cP)$ such that $\omega(\sigma)$ is not in $MP_0(N, \cP')$. By taking a decreasing sequence of neighborhoods whose intersection is $\rho$ we can find a sequence of representations $\sigma_i \rightarrow \rho$ such that the $\omega(\sigma_i)$ are not in $MP_0(N, \cP')$ for all $i$. However since $MP_0(N, \cP')$ is open in $AH(N, \cP')$ this contradicts Lemmas \ref{parabolics} and \ref{barcontinuous}. \qed{minimalparabolics}

We can now finish the definition of $\Psi$. Again there will be two cases depending on whether or not $\sigma(b)$ is parabolic. As before the case where $\sigma(b)$ is not parabolic will be where the work is.

Assume that $\sigma(b)$ is not parabolic. By Lemma \ref{standardconjugation} we can assume that
$$\omega(\sigma)(b) = \left(\begin{array}{cc}1 & 2 \\ 0 & 1 \end{array} \right).$$
Let $\beta$ be an essential loop in $\hat{M}_\sigma$ such that $\psi^{-1}_\sigma(\beta)$ bounds a disk in the Margulis tube $\Tu^{\epsilon_M}(\gamma_\sigma, M_\sigma)$. There will then be a unique element of $\pi_1(\hat{M}_\sigma)$ freely homotopic to $\beta$ whose $PSL_2\cx$ representative will be of the form
$$\left(\begin{array}{cc} 1 & z \\ 0 & 1 \end{array} \right)$$
with $\Im z>0$.
We then define
$$\Psi(\sigma) = \left\{ \begin{array}{cl} (\omega(\sigma), z) & \mbox{if $\sigma(b)$ is not parabolic} \\ (\omega(\sigma), \infty) & \mbox{if $\sigma(b)$ is parabolic.} \end{array} \right.$$ 

\subsection{Continuity, injectivity and inverses}
In this section we establish continuity properties of $\Phi$ and $\Psi$ and show that on suitably defined neighborhoods the maps are inverses of each other.

We begin with $\Phi$.
\begin{prop}
\label{phicontinuity}
The map $\Phi$ is continuous on $\Ao$.
\end{prop}

{\bf Proof.} First assume that $(\sigma, z)$ is a pair in $\Ao$ with $z \neq \infty$. Let $\widehat{\ab}$ be the Ahlfors-Bers parameterization for $MP(\hat{N}, \hat{\cP})$ and let $\ab$ be the Ahlfors-Bers parameterization for $MP(N, \cP)$. Note that the $\widehat{\ab}$-image of $\sigma_z = \tau(\sigma, z)$ and the $\ab$-image of $\Phi(\sigma,z)$ will be in isomorphic Teichm\"uller spaces and we can assume that these Teichm\"uller spaces have been identified such that $\Phi(\sigma,z) = \ab^{-1} \circ \widehat{\ab} \circ \tau(\sigma, z)$. Then the restriction of $\Phi$ to the component of $\int \cA_K$ that contains $(\sigma, z)$ is $\ab^{-1} \circ \widehat{\ab} \circ \tau$. From this continuity at $(\sigma, z)$ follows.

Now assume that $z = \infty$ and let $(\sigma_i, z_i)$ be a sequence with $\sigma_i \rightarrow \sigma$ and $z_i \rightarrow \infty$. Then the filling maps, $\phi_{\sigma_i, z_i}$, will have bi-Lipschitz constants limiting to $1$. As in the proof of Lemma \ref{barcontinuous} this will imply that $\Phi(\sigma_i, z_i)$ limits to $\Phi(\sigma, \infty) = \sigma$. \qed{phicontinuity}

We next address the continuity of $\Psi$.
\begin{prop}
\label{Psiinfinitycontinuity}
The map $\Psi$ is continuous all points in $\Vo \cap MP_0(N, \cP')$.
\end{prop}

{\bf Proof.} Let $\sigma_i$ be a sequence converging to $\sigma \in \Vo \cap MP_0(N, \cP')$ and let $z_i$ be the complex number in the pair $\Psi(\sigma_i)$. By Lemma \ref{barcontinuous} we already know that $\omega(\sigma_i) \rightarrow \sigma$. We need to show that $z_i \rightarrow \infty$. This follows from (2) of the drilling theorem since $\hat{M}_{\sigma_i}$ is the $\gamma_\sigma$-filling of $M_{\sigma_i}$, the length of the $\gamma_{\sigma_i}$ limits to zero and 
$$\frac{|z|}{\sqrt{2\Im z}}$$
is the normalized length of the meridian of $\hat{M}_\sigma$. \qed{Psiinfinitycontinuity}
%



Recall that we have fixed a neighborhood $V$ of $\rho$ upon which the holonomy inverse $h^{-1}$ is defined. Let $U'$ be the neighborhood of $g_\rho$ given by (2) of Theorem \ref{holonomyinverse}.
\begin{lemma}
\label{invertneighborhood}
There exists a neighborhood $U$ of $(\rho, \infty)$ in $\Ao$ with the following properties. Let $\Uo$ be the intersection of $U$ and $\Ao$. If $(\sigma,z)$ is a pair in $\Uo$ with $z \neq \infty$ then the pull back of the hyperbolic metric in $M_{\sigma,z}$ via $f_{\sigma, z}$ is in $U'$. If $z = \infty$ then the metric $g_\sigma$ is in $U'$.
\end{lemma}

{\bf Proof.} There exists a $K$ such that every metric in $\cH(N)$ that is $K$-bi-Lipschitz close to $g_\rho$ is in $U'$. We then choose our neighborhood $U$ such that for every $(\sigma, z) \in \Uo$ the metric $g_\sigma$ is $K/2$-bi-Lipschitz close to $g_\rho$ and such that $|z|$ is sufficiently large that the filling map $\phi_{\sigma, z}$ is $2$-bi-Lipschitz. \qed{invertneighborhood}

We now have a neighborhood where $\Phi$ is invertible.
\begin{prop}
\label{inverse}
The restriction of $\Psi \circ \Phi$ to $\Uo$ is the identity and therefore $\Phi$ is injective on $\Uo$.
\end{prop}

{\bf Proof.} Fix $(\sigma, z) \in \Uo$ and let $\bar{\sigma} = \Phi(\sigma, z)$. It follows directly from the definitions that $\Psi \circ \Phi(\sigma, \infty) = (\sigma, \infty)$ so we can assume that $z \neq \infty$.

We claim that $f_{\sigma, z}$ and $s_{\sigma_z}$ are homotopic as maps to  $M_{\sigma,z} - \gamma_{\sigma_z}$. By Lemma \ref{invertneighborhood} and (2) of Theorem \ref{holonomyinverse} there exists a submanifold $N_0$ of $N$, obtained by removing an open collar neighborhood of $\del N$, and a diffeomorphism
$$\alpha: N \longrightarrow N,$$
isotopic to the identity such that $s_{\sigma_z} \circ \alpha = f_{\sigma,z}$ on $N_0$. Therefore we have a homotopy of $s_{\sigma_z}$ to $f_{\sigma,z}$ where the domain of the homotopy is $N_0$. The image of the homotopy is contained in the union of $s_{\sigma_z}(N)$ and $f_{\sigma,z}(N)$ so $s_{\sigma_z}$ and $f_{\sigma,z}$ are homotopic as maps from $N_0$ to $M_{\sigma_z} - \gamma_{\sigma_z}$. Since the inclusion of $N_0$ in $N$ is a homotopy equivalence this implies that $s_{\sigma_z}$ are $f_{\sigma,z}$ are also homotopic as maps from $N$ to $M_{\sigma_z} - \gamma_{\sigma_z}$.


The proof is finished by directly applying the definition of $\Psi$. \qed{inverse}

We now have the following important corollary of Propositions \ref{Psiinfinitycontinuity} and \ref{inverse}.
\begin{cor}
\label{openlocalhomeo}
The $\Phi$ image of any neighborhood of $(\rho, \infty)$ in $\Ao$ contains a neighborhood of $\rho$ in $MP_0(N, \cP') \cup MP(N, \cP)$ and therefore $\Phi$ is a local homeomorphism from $\Ao$ to $MP_0(N, \cP') \cup MP(N, \cP)$ at $(\rho, \infty)$.
\end{cor}

\section{The punctured torus}

We now turn our attention to the punctured torus. We will begin by reviewing Minsky's fundamental work on the classification of Kleinian punctured torus groups. There is now a complete classification of all Kleinian surface groups but there are special features of this classification in the case of punctured torus groups that will be advantageous to us.

Let $\hat{T}$ be the punctured torus and $\bar{T}$ a compact torus with one boundary component. Define
$$N = \bar{T} \times [-1,1].$$
Let
$$\cP = \del \bar{T} \times [-1,1].$$
We fix a presentation of $\pi_1(\hat{T}) = \pi_1(N)$ by setting
$$\pi_1(N) = \langle a,b \rangle.$$
Of course, this is just the free group on 2 generators but $AH(N, \cP)$ is conjugacy classes of discrete faithful representations of the free group with the added condition that the image of the commutator $[a,b]$ is parabolic.

Our choice of generators for $\pi_1(\hat{T})$ also determines generators for $H_1(\hat{T})$. Using these generators we can write every element of $H_1(\hat{T})$ as a pair of integers. For each rational number $p/q \in \ratls \cup \infty$ let $\gamma_{p/q}$ be a simple closed curve that represents the unoriented homology class $\pm (q,p)$. In particular $\gamma_{0/1}$ is freely homotopic to $a$ and $\gamma_{1/0}$ is freely homotopic to $b$.

Let $\cT$ be the Teichm\"uller space of $\hat{T}$. Our choice of generators also canonically identifies $\cT$ with the upper half plane in $\cx$. As a subset of $\chat$ the upper half plane is compactified by $\reals \cup \infty$. This compactification is exactly the Thurston compactification by projective measured laminations. As above if $p/q \in \ratls \cup \infty$ then the lamination is the simple closed curve $\gamma_{p/q}$. If $\lambda \in \reals$ is not rational then the lamination is a minimal lamination without simple closed curves.



We denote the compactfied Teichm\"uller space by $\bar{\cT}$. Let $\Delta$ be the diagonal in $\del \bar{\cT} \times \del \bar{\cT}$. We then have the following important theorem of Minsky.
\begin{theorem}[Minsky \cite{Minsky:torus}]
\label{elc}
There exists a continuous bijection
$$\nu: \bar{\cT} \times \bar{\cT} - \Delta \longrightarrow AH(N, \cP)$$
and the composition $\nu \circ \ab$ is the identity on $MP(N, \cP)$.
\end{theorem}
Note that $MP(N, \cP)$ has only one component so $MP_0(N, \cP) = MP(N, \cP)$.

Although the map $\nu$ is not a homeomorphism if we restrict to certain subspaces it will become a homeomorphism. In particular the restriction of $\nu$ to $\cT \times \{pt\}$ gives an embedding of the Teichm\"uller space $\cT$ in $AH(N, \cP)$. If the point we choose is in the interior of $\cT$ then this is a {\em Bers' embedding} of Teichm\"uller space and the image is a {\em Bers' slice}. We can extend the map to the compactification $\bar{\cT}$ and $\nu$ restricted to $\bar{\cT} \times \{pt\}$ will be a homeomorphism because a continuous bijection between compact sets is always a homeomorphism.

Of more interest to us is the case when the point is a rational point, $\gamma_{p/q}$, in the boundary of $\bar{\cT}$. In this case the closure of $\cT \times \{\gamma_{p/q}\}$ in $\bar{\cT} \times \bar{\cT} - \Delta$ is not compact so the following result is not a direct corollary of Theorem \ref{elc}.

\begin{theorem}[Minsky \cite{Minsky:torus}]
\label{minskymaskit}
The restriction of $\nu$ to $(\bar{\cT} - \{\gamma_{p/q}\}) \times \{\gamma_{p/q}\}$ is a homeomorphism onto its image.
\end{theorem}

The image of $(\bar{\cT} - \{\gamma_{p/q}\}) \times \{\gamma_{p/q}\}$ has a natural embedding in $\cx$. We will restrict to the case when $p/q = 1/0$ but the construction clearly works for any $p/q$. Let $A$ be a closed annulus in $\bar{T}$ with core curve $\gamma_{1/0}$ and let $\cP'$ be the union of $A \times \{1\}$ and $\cP$. For each $z \in \cx$ we can define a representation $\sigma_z \in R(N, \cP')$ by setting
$$\sigma_z(a) = \left( \begin{array}{cc} \imath z & \imath \\ \imath & 0 \end{array} \right) \mbox{\ and \ } \sigma_z(b) = \left( \begin{array}{cc} 1 & 2 \\ 0 & 1 \end{array} \right).$$
If $\sigma_z$ is in $AH(N, \cP')$ we let $M_z$ be the hyperbolic manifold $\hthree/\sigma_z(\pi_1(N))$. By tameness (\cite{Bonahon:tame}) $M_z$ is homeomorphic to the interior of $N$.

We leave the following lemma as an exercise for the reader.
\begin{lemma}
\label{Nstandardform}
The map
$$z \mapsto \sigma_z$$
is a homeomorphism from $\cx$ to $R(N, \cP')$.
\end{lemma}

By Theorem \ref{elc} the space $MP(N, \cP')$ has exactly two components. Which of these components is $MP_0(N, \cP')$ is determined by the choice of orientation of $N$. We choose this orientation by picking a $\sigma_z$ in $MP(N, \cP')$ with $\Im z > 0$ and then fixing the orientation for $N$ such that $\sigma_z$ is in $MP_0(N, \cP')$.

Define the subset $\cM$ of $\cx$ to be the pre-image of $AH(N, \cP')$ of the map $z \mapsto \sigma_z$. We once again use Theorem \ref{elc} to see that $\cM$ has two components each canonically isomorphic to $\bar{\cT} -\gamma_{1/0}$. Let $\cM^+$ be the component of $\cM$ that contains the pre-image of $MP_0(N, \cP')$ and label the other component $\cM^-$.  The set $\cM^+$ (or its interior) is usually called the {\em Maskit slice}.

We now establish some elementary facts about the set $\cM$.
\begin{prop}
\label{Maskitfacts}
\begin{enumerate}
\item \label{symmetry} If $z$ is in $\cM^\pm$ then $z +2$ and $z-2$ are in $\cM^\pm$ and $-z$ and $\zbar$ are in $\cM^\mp$.


\item \label{notreal} The set $\cM$ does not intersect $\reals$ and therefore $\cM^+$ is contained in the upper half plane and $\cM^-$ is contained in the lower half plane.

\item \label{noline} The boundaries of $\cM^+$ and $\cM^-$ are not horizontal lines.
\end{enumerate}
\end{prop}

{\bf Proof (\ref{symmetry}).} The automorphism of $\pi_1(N)$ that sends $a \mapsto ba$ and $b \mapsto b$ sends $\sigma_z$ to $\sigma_{z+2}$ and therefore $\sigma_z$ is in $AH(N, \cP')$ if and only if $\sigma_{z+2} \in AH(N, \cP')$. This automorphism is induced by an orientation preserving diffeomorphism of $(N, \cP')$ so $z$ is in $\cM^\pm$ if and only if $z + 2$ is in $\cM^\pm$.

The automorphism of $\pi_1(N)$ that sends $a \mapsto a$ and $b \mapsto b^{-1}$ sends $\sigma_z$ to a representation conjugate to $\sigma_{-z}$. Note that this automorphism is induced by an orientation reversing diffeomorphism of $(N, \cP')$ so if $z$ is in $\cM^\pm$ then $-z$ is in $\cM^\mp$.

The automorphism that sends $a \mapsto a^{-1}$ and $b \mapsto b$ sends $\sigma_z$ to a representation conjugate to $\sigma_\zbar$. This is also induced by an orientation reversing diffeomorphism of $(N, \cP')$ so if $z \in \cM^\pm$ then $\zbar$ is in $\cM^{\mp}$.

{\bf (\ref{notreal}).}  If $x \in \reals$ then $\sigma_x$ is Fuchsian. That is $\sigma_x$ will fix $\reals$. The only Fuchsian groups in $AH(N, \cP)$ will lie in the image of the diagonal in $\cT \times \cT$ under $\nu$.

{\bf (\ref{noline}).} Given a representation $\sigma$ and an element $\gamma$ in $\pi_1(N)$ the trace of $\sigma(\gamma)$ is a complex number defined up to sign. This trace function is constant on conjugacy classes of $\pi_1(N)$ so it makes sense to discuss the trace of a curve $\gamma_{p/q}$. 
By Theorem \ref{minskymaskit} for each $p/q \in \ratls$ there is a unique $z \in \cM^+$ such that $\sigma_z(\gamma_{p/q})$ is parabolic and therefore its trace is $\pm 2$. Furthermore, this point will lie on the boundary of $\cM^+$. A similar statement holds for $\cM^-$.  We will find two rational values where the corresponding imaginary values of $z$ are different. This will finish the proof of the lemma.

For the curve $\gamma_{0/1}$ we need the trace of $\sigma_z(a)$ to be $\pm 2$ and therefore we must have $z = \pm 2\imath$. For the curve $\gamma_{1/2}$ we need trace of $\sigma_z(a^2b)$ to be $\pm 2$. This is equivalent to solving the quadratic polynomials
$$z(z+2) + 2 = \pm 2$$
so there are 4 solutions for $z$: $0$, $-2$ and $-1 \pm \imath \sqrt{3}$. As none of these solutions have $\Im z = \pm 2$ we are done. \qed{Maskitfacts}

{\bf Convention.} We let $\overline{\ab}$ be the Ahlfors-Bers parameterization of $MP(N, \cP')$. If we follow the definition of the Ahlfors-Bers parameterization given in Section \ref{abparam} then $\overline{\ab}$ is a map to two disjoint copies of the Teichm\"uller space $\cT$.  It will be convenient for us to identify these two copies and assume that $\overline{\ab}$ is a two-to-one map to $\cT$. By (\ref{symmetry}) we can also assume that $\overline{\ab}(\sigma_z) = \overline{\ab}(\sigma_\zbar)$.

\medskip

Given a $\sigma_z \in AH(N, \cP')$ we now make a few observations about the geometry of $M_z$ that will allows us to construct a marking map with some special properties. Since $\sigma_z(b)$ is parabolic the boundary of the convex core of $M_z$ will contain a triply-punctured sphere homeomorphic to $\hat{T} - \gamma_{1/0}$. The complement of the convex core will contain a component that meets this triply punctured sphere. Let $M'_z$ be the submanifold of $M_z$ obtained by removing this component of the complement of the convex core. Then $M'_z$ is a convex hyperbolic 3-manifold with boundary a triply punctured sphere. Next fix a hyperbolic metric on the surface $(\bar{T} - \int A) \times \{1\}$ so that the surface is isometric to the unique, finite area, complete, hyperbolic structure on the triply punctured sphere with the $\epsilon_M$-cusps removed.

We are now ready to define the marking map. Choose
$$f_z : N \longrightarrow M_z$$
so that $f_z$ restricted to $(\bar{T} - \int A) \times \{1\}$ is an isometry onto the boundary of the convex core. There are exactly six such isometries. Only one of these will give the correct marking of $M_z$ and we choose this one. By tameness we can assume that $f_z$ is an embedding whose image is contained the convex core of $M_z$ and that $f_z(N) \cap M_z^{\pr,\leq \epsilon} = f_z(\cP')$. Furthermore this map will be orientation preserving if $\Im z>0$ and orientation reversing if $\Im z<0$.


We now define $\hat{N}$ and $\hat{\cP}$ as in the previous section. In particular let $W$ be an open solid torus neighborhood of $\gamma_{1/0} \times \{0\}$  defined by $W = \int A \times (-1/2,1/2)$ and set $\hat{N} = N - W$. Then $\pi_1(\hat{N})$ has presentation
$$\langle a,b,c | bc = cb \rangle.$$
Let $\hat{\cP}$ be the union of $\cP$ and the boundary of $W$.

For each $(z,w) \in \cx \times \cx$ define a representation $\sigma_{z,w}$ of $\pi_1(\hat{N})$ by $\sigma_{z,w}(a) =\sigma_z(a)$, $\sigma_{z,w}(b) = \sigma_z(b)$ and
$$\sigma_{z,w}(c) = \left( \begin{array}{cc} 1 & w \\ 0 & 1 \end{array} \right).$$
Note $\sigma_{z,w} = \tau(\sigma_z, w)$ where $\tau$ is as defined in the previous section.

Again we let the reader check the following lemma.
\begin{lemma}
\label{standardform}
The map from $\cx \times \cx$ to $R(\hat{N}, \hat{\cP})$ defined by
$$(z,w) \mapsto \sigma_{z,w}$$
is a homeomorphism.
\end{lemma}
%

We now give a sufficient condition for this representation to be discrete and faithful.
\begin{lemma}
\label{sufficient}
The representation $\sigma_{z, w}$ is discrete and faithful if $z$ is in $\cM^+$ and $z-w$ is in $\cM^-$. Furthermore, if both $z$ and $z-w$ are in the interior of $\cM$ then $\sigma_{z,w}$ is in $MP_0(\hat{N}, \hat{\cP})$ and the Ahlfors-Bers parameterization, $\widehat{\ab}$, for $MP(\hat{N}, \hat{\cP})$ can be chosen such that $\widehat{\ab}(\sigma_{z,w}) = (\overline{\ab}(\sigma_z), \overline{\ab}(\sigma_{z-w}))$.
\end{lemma}

{\bf Proof.} We break $\hat{N}$ into two pieces. Let $N^+ \subset \hat{N}$ be those points in $\hat{N}$ where the second coordinate is non-negative. Then $N^-$ is defined similarly with non-negative replaced with non-positive. Let $\cP^+$ and $\cP^-$ be the restriction of $\hat{\cP}$ to  $N^+$ and $N^-$ respectively. Both $(N^-, \cP^-)$ and $(N^+, \cP^+)$ are homeomorphic to $(N, \cP')$. Choose a homeomorphism
$$f_+: (N^+, \cP^+) \longrightarrow (N, \cP')$$
and then define a homeomorphism
$$f_-:(N^-, \cP^-) \longrightarrow (N, \cP')$$
by $f_-(x, t) = f_+(x, -t)$.

Since $z$ and $z-w$ are both in $\cM$, the representations $\sigma_z$ and $\sigma_{z-w}$ are both in $AH(N, \cP')$. Then $M_z$ and $M_{z-w}$ are the corresponding quotient hyperbolic 3-manifolds and $M'_z$ and $M'_{z-w}$ the convex submanifolds bounded by totally geodesic triply-punctured spheres as defined above. Let $\hat{M}_{z,w}$ be the complete hyperbolic manifold obtained by gluing $M'_z$ and $M'_{z-w}$ along their totally geodesic boundaries. There are six possible such gluings. We choose the unique one such that $f_z \circ f_-$ and $f_{z-w} \circ f_+$ agree on the triply punctured sphere $\hat{T} - A \times \{0\} \subset \hat{N}$. We can then define a map
$$f_{z,w} : \hat{N} \longrightarrow \hat{M}_{z,w}$$
by
$$f_{z,w}(x,t) = \left\{ \begin{array}{cc} f_z \circ f_-(x,t) & \mbox{if $t \leq 0$} \\ f_{z-w} \circ f_+(x,t) & \mbox{if $t>0$.} \end{array} \right.$$

We now explicitly identify the group presentation
$$\langle a,b,c | [b,c] \rangle$$
with $\pi_1(\hat{N})$. Place the basepoint at a point in $\del W \cap \bar{T} \times \{-1/2\}$ and assume that $a$ and $b$ are the $\gamma_{0/1}$ and $\gamma_{1/0}$ curves on $\bar{T} \times \{-1/2\}$, respectively. We choose $c$ to bound a disk in $W$. To be more explicit let $c'$ be a proper arc in $A$ whose endpoints are in distinct components of $\del A$ and assume that the basepoint is contained in $c' \times \{-1/2\}$. Then $c = \del(c' \times [-1/2, 1/2])$. The final step is to choose orientations for the three curves. We will let the reader check that this can be done so that the marked hyperbolic 3-manifold $(\hat{M}_{z,w}, f_{z,w})$ induces the representation $\sigma_{z,w}$. The important point is the that the element $c^{-1} a$ in $\pi_1(\hat{N})$ will be freely homotopic to the curve $\gamma_{0/1}$ on $\hat{T} \times \{1/2\}$.

Note that the marking map we have constructed for $\sigma_{z,w}$ is an embedding and therefore $\sigma_{z,w}$ is in $AH_0(\hat{N}, \hat{\cP})$ (after choosing orientations appropriately). In particular if $z$ and $z-w$ are in the interior of $\cM$ then $\sigma_{z,w}$ is minimally parabolic and therefore is in $MP_0(\hat{N}, \hat{\cP})$. Furthermore, the components of the conformal boundary of $\hat{M}_{z,w}$ are exactly the punctured torus components of the conformal boundaries of $M_z$ and $M_{z-w}$ so we can define $\widehat{\ab}(\sigma_{z,w}) = (\overline{\ab}(\sigma_z), \overline{\ab}(\sigma_{z - w}))$ as desired.

\qed{sufficient}

{\bf Remark.} The proof of (1) of Proposition \ref{Maskitfacts} essentially showed that if $\overline{z-w} = z$ then there is an isometry from $M_z$ to $M_{z-w}$. This isometry can be used to construct an isometry
$$\iota: M_{z,w} \longrightarrow M_{z,w}$$
such that $\iota \circ f_{z,w} (x,t) = f_{z,w}(x,-t)$ from which it follows that $\iota_*(c) = c^{-1}$. We also note that $\overline{z-w} = z$ if and only if there is a $\mu \in \bar{\cT} - \gamma_{1/0}$ such that $\nu^{-1}(\sigma_z) = (\mu, \gamma_{1/0})$ and $\nu^{-1}(\sigma_{z-w}) = (\gamma_{1/0}, \mu)$. These two facts will play an important role in the proof of Lemma \ref{differentlaminations} below.

\medskip

The sufficiency conditions of the previous proposition can easily be adapted to be necessary conditions.

\begin{prop}
\label{necessary}
Given a pair $(z,w)$ with $\Im w \neq 0$ and let $s$ be the sign of $\Im w$.
The representation $\sigma_{z,w}$ is in $AH(\hat{N}, \hat{\cP})$ if and only if there exists an integer $n$ with
$$z- s n w \in \cM^+$$ and
$$z - s(n + 1)w \in \cM^-.$$
If $\sigma_{z,w} \in MP(\hat{N}, \hat{\cP})$ then we can choose the Ahlfors-Bers parameterization, $\widehat{\ab}$, such that $\widehat{\ab}(\sigma_{z,w}) = (\overline{\ab}(\sigma_{z-snw}), \overline{\ab}(\sigma_{z - s(n+1)w})$. The representation $\sigma_{z,w}$ is not in $AH(\hat{N}, \hat{\cP})$ if $\Im w = 0$. 
\end{prop}

{\bf Proof.} We will first show the necessity of the conditions by showing that if they don't hold there is a subgroup of $\pi_1(\hat{N})$ where the restriction of the representation is not discrete and faithful.

We first take care of the case where $\Im w = 0$. Then the restriction of $\sigma_{z,w}$ to the $\integers^2$ subgroup generated by $b$ and $c$ lies in $PSL_2\reals$. Since there are no discrete and faithful representations of $\integers^2$ in $PSL_2\reals$ the representation $\sigma_{z,w}$ cannot be discrete and faithful.

Now we assume that $\Im w \neq 0$. Note that for any integer $n$ the restriction of $\pi_1(\hat{N})$ to the subgroup generated by $\langle c^{-n} a, b \rangle$ is the representation $\sigma_{z-nw}$. If $z-nw$ is not in $\cM$ then $\sigma_{z-nw}$ is not discrete and faithful so neither is $\sigma_{z,w}$. Finding an $n$ such that $z-nw$ is not in $\cM$ will be our strategy for proving the necessity of the conditions.

There are two cases. If there exists an $n$ such that $\Im(z-nw) = 0$ then $z - nw$ is not in $\cM$ by (\ref{notreal}) of Proposition \ref{Maskitfacts}. If there is no such $n$ then there is a unique $n$ such that $\Im(z-snw) > 0$ and $\Im(z - s(n+1)w) < 0$. Now (\ref{notreal}) of Proposition \ref{Maskitfacts} implies that if  $z - snw$ is not in $\cM^+$ then it is not in $\cM$ and if $z - s(n+1)w$ is not in $\cM^-$ then it is not in $\cM$.


Now we prove sufficiency of the conditions. Assume $s = 1$ and define $z' = z - n w$.
Define an automorphism $\zeta$ of $\pi_1(\hat{N})$ by $\zeta(a) = c^{-n} a$, $\zeta(b) = b$ and $\zeta(c) = c$. Note that $\sigma_{z,w} \circ \zeta$ is the representation $\sigma_{z',w}$ which lies in $AH(\hat{N}, \hat{\cP})$ by  Lemma \ref{sufficient}. Therefore $\sigma_{z,w}$ is in $AH(\hat{N}, \hat{\cP})$. If $\sigma_{z,w}$ is in $MP(\hat{N}, \hat{\cP})$ then $\zeta$ takes the component $B$ of $MP(\hat{N}, \hat{\cP})$ that contains $\sigma_{z,w}$ to $MP_0(\hat{N}, \hat{\cP})$. We have defined $\widehat{\ab}$ for $MP_0(\hat{N}, \hat{\cP})$ in Lemma \ref{sufficient} and we use this to define $\widehat{\ab}$ on $B$ by 
$$\widehat{\ab}(\sigma_{z,w}) = \widehat{\ab}(\sigma_{z',w}) = (\overline{\ab}(\sigma_{z-nw}), \overline{\ab}(\sigma_{{z - (n+1)w}}).$$

If $s=-1$ the proof is exactly the same except that we define $\zeta(c) = c^{-1}$. \qed{necessary}


{\bf Remark.} Given a representation $\sigma_{z,w}$ in $AH(\hat{N}, \hat{\cP})$ we can find an immersion of $\hat{T}$ in $\hat{M}_{z,w}$ that induces the representation $\sigma_z$ of $\pi_1(\hat{T})$. This immersion will be homotopic to an embedding if and only if $n=0$ or $-1$ in the above lemma. In general the torus will wrap around the rank two cusp $n$-times if $n>0$ or $|n| -1$ times if $n <0$. It is the existence of such wrapped representations that leads to the bumping of components discovered in \cite{Anderson:Canary:pages}. This construction was also used in \cite{McMullen:graft} to show that quasifuchsian space self-bumped.

\medskip

We can now prove a version of Theorem \ref{elc} for $AH(\hat{N}, \hat{\cP})$. To do this we let $$\hat{\Delta} = (\{\gamma_{1/0}\} \times \bar{\cT}) \cup (\bar{\cT} \times \{\gamma_{1/0}\}).$$

\begin{theorem}
\label{elchat}
For each connected component $B$ of $AH(\hat{N}, \hat{\cP})$ there exists a homeomorphism
$$\hat{\nu}_B: \bar{\cT} \times \bar{\cT} - \hat{\Delta} \longrightarrow B$$
such that $\hat{\nu}_B \circ \widehat{\ab}$ is the identity.
\end{theorem}

{\bf Proof.} By Proposition \ref{necessary} for each $\sigma_{z,w}$ in $AH(\hat{N}, \hat{\cP})$ there will be a unique integer $n$ such that $z -snw$ is in $\cM^+$ and ${z - s(n+1)w}$ is in $\cM^-$ where $s$ is the sign of $\Im w$ (which is defined since $\Im w \neq 0$). Furthermore this $n$ and $s$ will be constant on  the entire connected component $B$.
We can define a homeomorphism from $B$ to $\cM^+ \times \cM^-$ by sending
$$\sigma_{z,w} \mapsto (z - snw, {z- s(n+1)w}).$$
The Maskit slice is canonically identified with $\hat{\cT} - \gamma_{1/0}$ so we can view the above map as a map to $\bar{\cT} \times \bar{\cT} - \hat{\Delta}$ and let $\hat{\nu}_B$ be the inverse of this map. Then Proposition \ref{necessary} implies that $\hat{\nu}_B \circ \widehat{\ab}$ is the identity on $\int B$ finishing the proof. \qed{elchat}

We now return to the maps $\Phi$ and $\Psi$ defined in the previous section. Recall that we have already defined $\Phi$ on the set $\Ao$ and in Proposition \ref{phicontinuity} we showed that $\Phi$ is continuous on $\Ao$. In the proof of Proposition \ref{phicontinuity} we saw that for points $(\sigma, z)$ in $\Ao$ with $z \neq \infty$ that $\Phi(\sigma, z)  = \ab^{-1} \circ \widehat{\ab} \circ \tau(\sigma, z)$. If $B$ is the component of $AH(\hat{N}, \hat{\cP})$ that contains $\tau(\sigma,z)$ then by Theorems \ref{elc} and \ref{elchat} we have $\Phi(\sigma, z) = \nu \circ \hat{\nu}_B^{-1} \circ \tau(\sigma,z)$. This observation is the key to extending $\Phi$ to all of $\cA_K$.

The first step is to show that the image of $\cA_K$ under the map $\hat{\nu}_B^{-1} \circ \tau$ lies in the domain of $\nu$.
\begin{lemma}
\label{differentlaminations}
Let $B$ be a component of $AH(\hat{N}, \hat{\cP})$ and let $C = \tau^{-1}(B) \cap \cA_K$. Then if $(\sigma, z) \in C$ the point $\hat{\nu}_B^{-1} \circ \tau(\sigma, z)$ is not in $\Delta$.
\end{lemma}

{\bf Proof.} Let $\hat{M}$ be the hyperbolic 3-manifold corresponding to the representation $\tau(\sigma,z)$ and let $\beta$ be a simple closed curve on the torus component of $\hat{\cP}$ that is freely homotopic to the element $c$ in $\pi_1(\hat{N})$. By the remark following the proof of Lemma \ref{sufficient} we see that if $\hat{\nu}_B^{-1} \circ \tau(\sigma, z)$ is in $\Delta$ then there is an isometry of $\hat{M}$ to itself and this isometry won't change the free homotopy class of $\beta$. 
We also know that the normalized length of $\beta$ is greater than $K$ (which we can assume is greater than $2\pi$) so by the Gromov-Thurston $2\pi$-Theorem (see \cite{Bleiler:Hodgson:twopi}) we can $\beta$-fill $\hat{M}$ to obtain a manifold $M$ with a a metric of negative sectional curvature that is equal to the original hyperbolic metric outside of the Margulis tube for the rank two cusp. The Gromov-Thurston construction is symmetric so the manifold $M$ will also have an isometry flipping the two ends of the manifold. This implies that the two (relative) ends of $M$ will have the same ending lamination which is a contradiction (\cite{Bonahon:tame}, \cite{Thurston:book:GTTM}).

We can actually see this more concretely in our particular situation. As noted $M$ is homeomorphic $\hat{T} \times (-1,1)$ and we can assume that the product structure is chosen such that the self-homeomorphisim of $M$ defined by $(x,t) \mapsto (x,-t)$ is an isometry. Since this map is homotopic to the identity this would imply that all closed geodesics in $M$ are contained in $\hat{T} \times \{0\}$ which is a contradiction. \qed{differentlaminations}

We can now extend $\Phi$ to all of $C$.
\begin{lemma}
\label{componentext}
Let $B$ be a component of $AH(\hat{N}, \hat{\cP})$ and let $C = \tau^{-1}(B) \cap \cA_K$. Then $\nu \circ \nu_B^{-1} \circ \tau$ is defined and continuous on all of $C$ and is equal to $\Phi$ on the interior of $C$.
\end{lemma}

{\bf Proof.} By Lemma \ref{differentlaminations} the $\hat{\nu}_B^{-1} \circ \tau$-image of $C$ is contained in the domain of $\nu$ and therefore $\nu \circ \hat{\nu}_B^{-1} \circ \tau$ is defined on all of $C$. As all three maps are continuous so is their composition. The interior of $C$ is mapped to the interior of $B$ by $\tau$. In particular the $\tau$-image of the interior of $C$ is contained in a component of $MP(\hat{N}, \hat{\cP})$ and therefore on the interior of $C$ we have
$$\nu \circ \hat{\nu}_B^{-1} \circ \tau = \ab^{-1} \circ \widehat{\ab} \circ \tau = \Phi.$$
\qed{componentext}
%

Applying Lemma \ref{componentext} to each component of $AH(\hat{N}, \hat{\cP})$ allows us to extend $\Phi$ to all of $\cA_K$.
\begin{cor}
\label{completecontinuityPhi}
The maps $\Phi$ extends continuously to all of $\cA_K$.
\end{cor}

Recall that we have defined a neighborhood $U$ of $(\rho, \infty)$ such that $\Phi$ is injective on $\Uo = U \cap \Ao$. We now see that $\Phi$ is injective on all of $U$.
\begin{prop}
\label{completeinjectivity}
The map $\Phi$ is injective on $U$.
\end{prop}

{\bf Proof.} We begin by observing that if $\sigma$ is in $MP_0(N, \cP')$ then $\Phi(\sigma',z) = \sigma$ if and only if $\sigma = \sigma'$ and $z = \infty$. If $z = \infty$ then this follows directly from the definition of $\Phi$. On the other hand if $z \neq \infty$ then $\hat{\nu}^{-1} \circ \tau(\sigma,z)$ is not in $\hat{\Delta}$ and therefore $\Phi(\sigma,z)$ is not in $MP_0(N, \cP')$.

We now check for injectivity at points $(\sigma, z)$ with $z \neq \infty$. As in the proof of Proposition \ref{completecontinuityPhi} we will use the fact that at such points the map $\Phi$ factors through a map to $\bar{\cT} \times \bar{\cT}$.  The second half of this factorization is the map $\nu$ which is injective so to prove the theorem we need to show that the map to $\bar{\cT} \times \bar{\cT}$ is injective.

Let $(\sigma, z)$ be a pair in $U$ and let $B$ be the component of $AH(\hat{N}, \hat{\cP})$ that contains $\tau(\sigma,z)$. Let $\mu = \hat{\nu}_B^{-1} \circ \tau(\sigma, z)$ and choose a sequence $\mu_i$ in $\cT \times \cT$ that converges to $\mu$. Since $\tau$ is a local homeomorphism at $(\sigma, z)$, for large $i$ we can find $\sigma_i$ and $z_i$ such that the pairs $(\sigma_i, z_i)$ are contained in $\Uo$, converge to $(\sigma,z)$, the $\tau(\sigma_i, z_i)$ are contained in $B$ and  $\hat{\nu}_B^{-1} \circ \tau(\sigma_i, z_i) = \mu_i$.

Let $(\sigma', z')$ be another pair in $U$ such that $\hat{\nu}_{B'}^{-1} \circ \tau(\sigma', z') = \mu$ where $B'$ is the component of $AH(\hat{N}, \hat{\cP})$ that contains $\tau(\sigma', z')$. By repeating the construction we can find pairs $(\sigma'_i, z'_i)$ in $\Uo$ as above. In particular, $\hat{\nu}_{B'}^{-1} \circ \tau (\sigma'_i, z'_i) = \mu_i = \hat{\nu}_B^{-1} \circ \tau(\sigma_i, z_i)$ which implies that $\Phi(\sigma'_i, z'_i) = \Phi(\sigma_i, z_i)$. Since $\Phi$ is injective on $\Uo$ it follows that $(\sigma'_i, z'_i) = (\sigma_i, z_i)$. Since the two sequences are identical they must have the same limit so $(\sigma', z') = (\sigma, z)$ and $\Phi$ is injective on $U$. \qed{completeinjectivity}

We now show that $\cA$ and $AH(N, \cP)$ are locally homeomorphic.
\begin{theorem}
\label{Philocalhomeo}
The map $\Phi$ is a local homeomorphism from $\cA_K$ to $AH(N, \cP)$ at $(\rho, \infty)$.
\end{theorem}

{\bf Proof.} Choose a compact neighborhood $C$ of $(\rho, \infty)$ that is contained in $U$. Since $\Phi$ is continuous and injective on $C$ the restriction of $\Phi$ to $C$ is a homeomorphism onto its image. Let $\Co =  C \cap \Ao$. By Corollary \ref{openlocalhomeo}, $\Phi(\Co)$ contains a neighborhood of $\rho$ in $MP_0(N, \cP') \cup MP(N, \cP)$. Note that $MP(N, \cP)$ is dense in $AH(N, \cP)$ so $\Phi(C) = \overline{\Phi\left(\Co\right)}$ contains a neighborhood of $\rho$ in $AH(N, \cP)$ and therefore $\Phi$ is a local homeomorphism at $(\rho, \infty)$. \qed{Philocalhomeo}

To finish the proof that $AH(N, \cP)$ is not locally connected we need to find a pair $(\rho, \infty)$ where $\cA$ is not locally connected. This is a fairly easy consequence of the facts about the Maskit slice given by Proposition \ref{Maskitfacts} and the description of $\cA$ given by Proposition \ref{necessary}.

To clarify the situation we expand on the comments made at the end of the introduction by stating a few general facts about subsets of $\cx$ and $\chat$ and then observe that they apply to our specific setting. First let $S$ be a subset of $\chat$ that is translation invariant (if $z \in S$ then $z \pm 2 \in S$) and contains the point $\infty$.  Then $S$ will not be locally connected at $\infty$ if the restriction of $S$ to $\cx$ contains a bounded component (in $\cx$).

For our second general observation we let $S^+$ and $S^-$ be translation invariant subsets of $\cx$ such that $S^+$ is contained in the upper half plane and $S^-$ is contained in the lower half plane. We also assume that these sets are not bounded by horizontal lines. Then we can translate $S^-$ such that the intersection of $S^+$ with the translation of $S^-$ is a translation invariant subset that contains bounded components.

Combining these two observations it is easy to find $z \in \cM^+$ such that $\cA_{\sigma_z}$ is not locally connected at $\infty$. This is not quite enough to see that $\cA$ is not locally connected. For this we need a somewhat stronger statement which we record in the following lemma.

\begin{lemma}
\label{rectangle}
There exists an open subset $\cO$ of $\cM^+$ and a closed rectangle $R$ such for every $z \in \cO$ the set $\cA_{\sigma_z}$ contains points in the interior of $R$ but $\cA_{\sigma_z} \cap \del R  = \emptyset$. We can choose $R$ such that it sides are parallel to the axes and its width is $<2$.
\end{lemma}

{\bf Proof.} Using \eqref{noline} of Proposition \ref{Maskitfacts} we can find a rectangle $Q$ and a point $z$ in the interior of $\cM^+$ with the following properties:
\begin{enumerate}
\item The sides of $Q$ are parallel to the axes and $Q$ is contained in the upper half plane. The two vertical sides and the lower horizontal side are disjoint from $\cM^+$. The width of $Q$ is $<2$.

\item The point $z$ is contained in the interior of $Q$ and the distance from $z$ to the top horizontal boundary component of $Q$ is exactly twice the distance to the lower horizontal boundary component of $Q$.
\end{enumerate}

Now let $R$ be the rectangle
$$R = \{w \in \cx | 3z - w \in Q\}.$$
One can check that for  $w \in R$, $\Im(3z-w)>0$ and $\Im(3z-2w)<0$ and therefore by Proposition \ref{necessary} a point $w \in R$ is in $\cA_{\sigma_{3z}}$ if and only if $3z - w \in \cM^+$ and $3z - 2w \in \cM^-$. We also note that by \eqref{symmetry} of Proposition \ref{Maskitfacts}, $3z-2w \in \cM^-$ if and only if $2w - 3z \in \cM^+$. It is now straightforward to check that $2z$ is in the interior of $R$ and in $\cA_{\sigma_{3z}}$ but no point in $\del R$ is in $\cA_{\sigma_{3z}}$. For the latter statement we note that for points $w$ in the upper horizontal or vertical sides of $\del R$, the point $3z-w$ is in the lower horizontal or vertical sides of $\del Q$. For a point $w$ in the lower horizontal side of $\del R$, the point $2w - 3z$ is in the lower horizontal side of $\del Q$.

To finish the proof we note that there is an $\epsilon$-neighborhood of the lower horizontal and vertical sides of $\del Q$ that is disjoint from $\cA_{\sigma_{3z}}$. We choose $\cO$ to be an open neighborhood of $3z$ in $\cM^+$ that is contained in the $\epsilon$-neighborhood of $3z$ in $\cx$. \qed{rectangle}

We can now prove the main theorem of the paper.

\begin{theorem}
\label{Anotlocalatmu}
There exists $\sigma \in MP_0(N, \cP')$ such that $AH(N, \cP)$ is not locally connected at $\sigma$.
\end{theorem}

{\bf Proof.} By Theorem \ref{Philocalhomeo} if we find a $z \in \cM^+$ such that $\cA$ is not locally connected at $(\sigma_z, \infty)$ then $AH(N, \cP)$ is not locally connected at $\sigma_z = \sigma$.  Let $R$ be the rectangle and $\cO$ the open subset of $\cM^+$ given by Lemma \ref{rectangle}. We claim that for all $z \in \cO$ the set $\cA$ is not locally connected at $(\sigma_z, \infty)$.

We identify $\cO$ with a subset $\cO_\sigma$ of $MP_0(N, \cP')$ via the map $z \mapsto \sigma_z$. Define an open neighborhood $U$ of $(\sigma_z, \infty)$ by taking the intersection of $\cO \times \chat$ with $\cA$. We will show that any other neighborhood $V \subset U$ of $(\sigma_z, \infty)$ has an infinite number of components.

Since $\cA$ has the subspace topology inherited from $MP_0(N, \cP') \times \chat$, the open set $V$ is the intersection with $\cA$ of an open set $V'$ in $\cO \times \chat$. The restriction of $V'$ to the plane $\sigma_z \times \chat$ will contain a neighborhood of infinity and hence an infinite number of translates of $R$ under the group generated by $z \mapsto z+2$. Each of these disjoint rectangles will contain a point in $\cA_{\sigma_z}$ and hence in $V$. Let $w_0$ and $w_1$ be two such points contained in distinct translations $R_0$ and $R_1$ of $R$. The component of $V$ containing $w_i$ will be contained $R_i \times \cO$. Since $R_0 \times \cO$ is disjoint from $R_1 \times \cO$, the points $w_0$ and $w_1$ are in distinct components of $V$. Therefore $V$ has an infinite number of connected components and $\cA$ is not locally connected at $(\sigma_z, \infty)$. \qed{Anotlocalatmu}



{\bf Remark.} We can also use Theorem \ref{Philocalhomeo} to find representations where $AH(N, \cP)$ is locally connected by finding pairs in $(\sigma_z, \infty)$ in $\cA$ where $\cA$ is locally connected. Here is one way to do this. Choose a $z$ sufficiently close to the boundary of $\cM^+$ such that for any $w \in \cx$ with $\Im w > 0$ then $z - nw \in \cM^+$ and $z -(n+1)w \in \cM^-$ implies that $n=0$. If this is the case then Proposition \ref{necessary} implies that $\cA_{\sigma_z}$ will just be a translated copy of $\cM^+$ and it is fairly easy to see that $\cA$ is locally connected at $(\sigma_z, \infty)$. In fact $\cA$ will not self-bump at $(\sigma_z, \infty)$. In the terminology of \cite{Holt:Souto:torus} the above condition is equivalent to the representation $\sigma_z$ not {\em wraping} (see the remarks after Proposition \ref{necessary}) and the fact that $AH(N, \cP)$ does not self-bump is a special case of the main theorem of \cite{Holt:Souto:torus}. We further note that by working a bit harder one can find representations where there is self-bumping but $AH(N, \cP)$ is still locally connected.

\section{Other deformation spaces}
In this brief final section we make some conjectures about the behavior of other deformation spaces.

The most obvious next case to consider is surfaces of higher genus.
\begin{conj}
\label{highergenus}
Let $S$ be a compact surface with boundary. Then $AH(S \times [0,1], \del S \times [0,1])$ is not locally connected.
\end{conj}
In fact one might make the stronger conjecture that for any pared 3-manifold that has an essential cylinder the deformation space is not locally connected.

Another natural deformation space to look at is the Bers' slice.
\begin{conj}
\label{bersslice}
If $S$ is not a punctured torus or a 4-times punctured sphere then any Bers' slice for $S$ is not locally connected.
\end{conj}
A Bers' slice is very similar to the deformation space of an acylindrical manifold and one could make a similiar conjecture for such manifolds. Of these two conjectures the second is probably more difficult and less likely to be true.

\bibliographystyle{math}
\bibliography{math}

\begin{sc}
\noindent
Department of Mathematics\\
University of Utah\\
155 S 1400 E, JWB 233\\
Salt Lake City, UT 84112
\end{sc}

 \end{document}